\documentstyle[12pt]{article}

\newcommand{\sgroup}[2]{\raisebox{-2.2ex}
{\begin{picture}(28,24)(0,0)
\put(4,15){\vector(1,0){20}}
\put(4,10){\vector(1,0){20}}
\put(10,3){${\scriptstyle #1}$}
\put(10,18){${\scriptstyle #2}$}
\end{picture}}}

\newcounter{th}

\newtheorem{th}[th]{Th\'{e}or\`{e}me}
\newtheorem{prop}[th]{Proposition}
\newtheorem{cor}[th]{Corollaire}
\newtheorem{lem}[th]{Lemme}
\newtheorem{rem}[th]{Remarque}

\newtheorem{exem}[th]{Exemple}

\newcounter{ith}

\newtheorem{ith}[ith]{Th\'{e}or\`{e}me}

\newtheorem{iexems}[ith]{Exemples}

\title{PR\'EQUANTIFICATION DE CERTAINES
VARI\'ET\'ES DE POISSON}
\author{F. Alcalde Cuesta\thanks{Recherche partiellement support\'ee par
D.G.I.C.Y.T. Espagne (Proyecto PB90-0765) et Xunta de Galicia (Proxecto
XUGA20704B90)} \\
{\normalsize Dpto. de Xeometria e Topoloxia, Universidade de Santiago} \\
{\normalsize 15706 Santiago de Compostela (Espagne)}}
\date{}
\begin{document}
\maketitle

\begin{abstract}
A surjective submersion $\pi : M \to B$ carrying a field of simplectic
structures on the fibres is {\em symplectic}
if this Poisson structure is {\em minimal}.
A symplectic submersion may be interpreted as a family of
mechanical systems depending on a parameter in $B$.
We give some conditions to find a closed form which represent the
foliated form $\sigma$ gluing the symplectic forms on the fibres.
This is the first step to prequantize all these systems at once.
We will indeed exhibit an integrality condition which does not
depend on the closed form representing $\sigma$: if the fibres are
1-connected and $H^3(B;\mbox{\bf Z})=0$, then
there exists a $S^1$-principal fibre bundle with a connection
whose curvature represents $\sigma$ iff the {\em group of
spherical periods of} $\sigma$ is a discrete subgroup of {\bf R}.

The {\it symplectic integration} of a Poisson manifold $(M,\Lambda)$
is a symplectic groupoid $(\Gamma,\eta)$ with 1-connected fibres such that the
space of units
with the induced Poisson structure is isomorphic to $(M,\Lambda)$.
This notion was introduced by A. Weinstein in order to quantize
Poisson manifolds by quantizing their symplectic integration.
We show that if the symplectic integration is
prequantizable, then there exists a unique prequantization which is trivial
over $M$.
Such a prequantization is a central extension of $\Gamma$ by $S^1$.
We show that the symplectic integration of a minimal Poisson manifold
is prequantizable iff
the group of spherical periods is discrete. Moreover we prove that a
{\em totally aspherical} Poisson manifold (any vanishing cycle is trivial and
the $\pi_2$ of the leaves is zero) is prequantized in the sense of Weinstein
by a trivial fibre bundle.
\end{abstract}

\section{Introduction}

\hspace{5mm}
Un fibr\'e localement trivial
\( \pi : M \rightarrow B \)
est {\em symplectique} \cite{GLSW} si le groupe des diff\'eomorphismes
de la fibre $F$ admet une r\'eduction au groupe des symplectomorphismes d'une
structure symplectique $\sigma_F$ sur $F$.
Un tel fibr\'e peut \^etre interpr\'et\'e
comme une famille de syst\`emes m\'ecaniques
qui d\'epend d'un param\`etre dans $B$. Il est naturel de s'int\'eresser \`a
la pr\'equantification globale de ces syst\`emes. On cherche \`a
construire un $S^1$-fibr\'e principal $q:E \to M$ muni d'une connexion $\theta$
dont la courbure
$\omega$ soit une extension ferm\'ee des formes symplectiques des fibres.
Dans \cite{GLSW}, on exhibe un crit\`ere d'existence de telles extensions et
l'on d\'emontre que  ce crit\`ere est effectivement v\'erifi\'e si $B$ et  $F$
sont  connexes et simplement connexes.  Dans ces conditions, pour qu'il existe
une pr\'equantification globale, il faut et il suffit que le  groupe des
p\'eriodes de $\sigma_F$ soit un sous-groupe discret de {\bf R}.  C'est la
condition  d'int\'egrabilit\'e qui caract\'erise l'existence d'une
pr\'equantification  de $F$ (cf. \cite{So} et \cite{TW}).

Les fibres de $\pi$ d\'efinissent un feuilletage $\cal F$ et
leurs formes symplectiques se re\-collent en
une {\em forme feuillet\'ee symplectique} $\sigma$. Le couple
$({\cal F},\sigma)$
d\'etermine donc une structure de Poisson $\Lambda$ sur $M$.
Dans \S\ref{Smin}, on montrera que $\Lambda$ est
{\em minimale} au sens de \cite{DH}. D'autre part,
les obstructions de \cite{DH} et \cite{V1} \`a
l'existence d'un repr\'esentant ferm\'e de $\sigma$ permettront
de reformuler la condition de \cite{GLSW} (voir \S\ref{obst}).
\vspace{2ex}

Cette remarque situe la pr\'equantification des fibr\'es symplectiques dans le
cadre g\'en\'eral de la quantification g\'eom\'etrique des vari\'et\'es de
Poisson. On connait des diff\'erents proc\'ed\'es:
\vspace{1ex}

\noindent
1) celui de I. Vaisman dans \cite{V2} en termes de fibr\'es en droites
complexes et d\'erivations contravariantes;
\vspace{1ex}

\noindent
2) le programme de quantification de A. Weinstein (\cite{W}): on r\'ealise
la vari\'et\'e de Poisson $(M,\Lambda)$ comme l'espace des unit\'es d'un
groupo\"{\i}de symplectique $(\Gamma,\eta)$  \`a fibres connexes et simplement
connexes,   que l'on appelle {\em int\'egration symplectique}, puis on
quantifie
$(M,\Lambda)$ en quantifiant $(\Gamma,\eta)$;
\vspace{1ex}

\noindent
3) la version infinit\'esimale de J. Huebschmann dans \cite{H}.
\vspace{2ex}

La pr\'equantification des fibr\'es symplectiques sugg\`ere une approche
directe dans la premi\`ere \'etape de la quantification des vari\'et\'es de
Poisson: le recollement des  pr\'equantifications des feuilles symplectiques.
Soit $\Lambda=({\cal F},\sigma)$ une structure de Poisson r\'eguli\`ere sur une
vari\'et\'e $M$.  Une {\em pr\'equantification} de $(M,\Lambda)$ est
la donn\'ee d'un $S^1$-fibr\'e principal
\( q : E \rightarrow M \)
et d'une connexion $\theta$ dont la courbure $\omega$ repr\'esente la
forme feuillet\'ee symplectique $\sigma$. Les obstructions de \cite{DH}
et \cite{V1} sont encore des obstructions \`a la pr\'equantification. En
particulier, la vari\'et\'e de Poisson  $(M,\Lambda)$ doit \^etre minimale. Par
ailleurs, si l'on fixe un repr\'esentant ferm\'e $\omega$ de  $\sigma$,
l'existence d'une pr\'equantification  est caract\'eris\'ee par les p\'eriodes
de $\omega$. Mais, m\^eme dans le cas  des fibr\'es symplectiques, la condition
d'int\'egrabilit\'e ne restera pas valable  si l'on change de repr\'esentant
ferm\'e.  On est donc amen\'e \`a d\'efinir et utiliser des p\'eriodes
(sph\'eriques) de $\sigma$ qui  ne d\'ependent pas du choix de repr\'esentant
ferm\'e.

Toute sph\`{e}re point\'ee contenue dans une feuille de $\cal F$ peut \^{e}tre
pouss\'ee dans le feuilles voisines. Par int\'{e}gration de $\sigma$
sur ces sph\`{e}res, on obtient une {\em fonction d'aire} sur une transversale.
Ces fonctions d'aire d\'efinissent  un sous-groupo\"{\i}de $Per(\sigma)$ de $M
\times \mbox{\bf R}$ appel\'e le  {\em groupo\"{\i}de des p\'eriodes
sph\'eriques} \cite{AH} de $\sigma$.  Si l'on oublie les points base, ces
p\'eriodes engendrent un sous-groupe $\cal P$ de {\bf R} que  l'on appelera le
{\em groupe des  p\'eriodes sph\'eriques de $\sigma$}.  Ce groupe jouera le
r\^ole du groupe des p\'eriodes dans la condition d'integrabilit\'e. Mais la
construction d'une pr\'equantification exige aussi l'annulation  des
obstruction
cohomologiques: c'est une condition n\'ecessaire pour le recollement des
pr\'equantifications des feuilles. Pour cela, il faut disposer au pr\'ealable
de
r\'esultats  cohomologiques analogues \`a ceux de \cite{GLSW} pour les
fibr\'es.  \vspace{2ex}

Les r\'esultats de \cite{AH} sur la structure cohomologique des submersions
surjectives (que l'on rappelera au \S\ref{annul}) vont permettre d'\'et\'endre
les r\'esultats de \cite{GLSW}. Ainsi, on  d\'emontrera au \S\ref{repfer} un
crit\`ere d'existence de repr\'esentants  ferm\'es. Une submersion  surjective
\( \pi : M \rightarrow B \)  munie d'une forme feuillet\'ee symplectique
$\sigma$ est  {\em symplectique} si la structure de Poisson correspondante est
minimale.  Dans ce contexte, on obtiendra la condition de pr\'equantification
suivante:

\begin{ith}    \label{condint}
Soit
\( \pi : M \rightarrow B \)
une submersion symplectique dont les fibres sont connexes et simplement
connexes et$H^3(B;\mbox{\bf Z})=0$. Alors, la submersion symplectique est
pr\'equantifiable si, et seulement si,  $\cal P$ est un sous-groupe discret de
{\bf R}.
\end{ith}

Les conditions du th\'eor\`eme~\ref{condint} entra\^{\i}nent
l'annulation des obstructions
(voir \S\ref{sub}). Si $\cal P$ est de plus un sous-groupe discret de {\bf R},
il est possible d'int\`egrer $\sigma$ en un cocycle sur $M$ \`a valeurs dans
$\mbox{\bf R}/{\cal P}$
(voir \S\ref{th1}). En fait, ce proc\'ed\'e d'int\'egration cohomologique
fournira un {\em bon} repr\'esentant ferm\'e de $\sigma$ dont le groupe des
p\'eriodes est \'egal \`a $\cal P$.  \vspace{2ex}

Dans le cas g\'en\'eral des feuilletages, les r\'esultats cohomologiques de
\cite{AH} ne restent pas valables: la cohomologie d'un feuilletage en droites
irrationnelles du tore $T^2$ n'est pas triviale  (voir \cite{He}), bien que les
feuilles sont contractiles. Donc la preuve du  th\'eor\`eme~\ref{condint}
pousse
jusqu'au bout la d\'emarche de pr\'equantification de \cite{GLSW}.
N\'eanmoins,
le passage \`a l'int\'egration symplectique ram\`enera l'\'etude cohomologique
(et donc le probl\`eme de la pr\'equantification) des vari\'et\'es  de Poisson
r\'eguli\`eres au cas des submersions \`a fibres connexes et simplement
connexes.  D'autre part, l'emploi de l'int\'egration symplectique permet de
reformuler le probl\`eme de la  pr\'equantification: on cherche \`a construire
une extension centrale d'un groupo\"{\i}de  symplectique (cf. \cite{WX}). Il
s'agit de g\'en\'eraliser les  r\'esultats de \cite{TW} sur les extensions
centrales de groupes.

Soit $(M,\Lambda)$ une vari\'et\'e de Poisson int\'egrable.
Une {\em pr\'equantification au sens de Weinstein} de $(M,\Lambda)$ est une
pr\'equantification $ q : (E,\theta) \rightarrow (\Gamma,\eta)$ de son
int\'egration symplectique. Puisque  $M$ est une sous-vari\'et\'e lagrangienne
de $\Gamma$, ce fibr\'e induit  un fibr\'e plat au-dessus de $M$. Si
l'holonomie
de la connexion  induite est triviale, celle-ci trivialise le fibr\'e induit et
l'on dira que la  pr\'equantification est {\em triviale en restriction \`a
$M$}.   Dans \S\ref{extg}, on prouvera le th\'eor\`eme suivant:

\begin{ith}              \label{canon}
Soit $(M,\Lambda)$ une vari\'et\'e de Poisson pr\'equantifiable au sens de
Weinstein. Alors, il existe
une unique pr\'equantification $ q : (E,\theta) \to (\Gamma,\eta) $ triviale
en restriction \`a l'espace des unit\'es $M$. En outre, l'espace total $E$ est
muni d'une structure canonique de groupo\"{\i}de de Lie
qui en fait une extension de $\Gamma$ par $S^1$.
\end{ith}

Le th\'eor\`eme~\ref{canon} a \'et\'e d\'emontr\'ee
dans \cite{WX} pour les groupo\"{\i}des symplectiques localement triviaux.
\vspace{2ex}

D'apr\`es \cite{DH}, une vari\'et\'e de Poisson minimale est toujours
int\'egrable. Dans \S\ref{th3}, on obtiendra la condition d'int\'egrabilit\'e
suivante:

\begin{ith}    \label{preW}
Soit $\Lambda=({\cal F},\sigma)$ une structure de Poisson minimale sur
une vari\'et\'e $M$ dont tout cycle \'evanouissant de $\cal F$ est
trivial. La vari\'et\'e de Poisson $(M,\Lambda)$ est pr\'equantifiable
au sens de Weinstein si, et seulement si, le groupe $\cal P$ des p\'eriodes
sph\'eriques de $\sigma$ est un sous-groupe discret de {\bf R}.
\end{ith}

Une structure de Poisson r\'eguli\`ere $\Lambda = ({\cal F},\sigma)$ est
{\em totalement asph\'erique} \cite{DH} si tout cycle \'evanouissant de
$\cal F$ est trivial et le $\pi_2$ des feuilles est nul.
Dans ce cas, le groupe $\cal P$ est nul et l'on aura le th\'eor\`eme
suivant:

\begin{ith}   \label{tot}
Si $(M,\Lambda)$ est une vari\'et\'e de Poisson totalement asph\'erique, alors
$(M,\Lambda)$ est pr\'equanti\-fi\'ee au sens de Weinstein par un fibr\'e
trivial en cercles.
\end{ith}

\begin{iexems} \label{exei}
{\em
Soit $\Lambda=({\cal F},\sigma)$ une structure de Poisson
sur une vari\'et\'e compacte $M$ de dimension 3. Si $\cal F$ est
transversalement orientable,
d'apr\`es le th\'eor\`eme de stabilit\'e de Reeb et le th\'eor\`eme de Novikov,
il y a trois cas possibles:

\noindent
1) $\cal F$ est totalement asph\'erique et donc $(M,\Lambda)$ est
pr\'equantifiable au sens de Weinstein d'apr\`es le th\'eor\`eme~\ref{tot};

\noindent
2) $\cal F$ est la fibration triviale en sph\`eres. Pour que cette fibration
soit symplectique, il faut et il suffit que la fonction d'aire
$p = - \!\!\!\!\! \int_{S^2} \sigma$ soit constante.
Une telle fibration est toujours pr\'equantifi\'ee par un fibr\'e principal
de groupe
$\mbox{\bf R}/ p \mbox{\bf Z}$. L'espace total est \'egal \`a $S^3 \times S^1$
et la longeur des  fibres est \'egale \`a $p$. D'habitude, on cherche \`a
pr\'equantifier par un  fibr\'e principal de groupe $\mbox{\bf R}/2 \pi \hbar
\mbox{\bf Z}$, o\`u $2 \pi \hbar$ est  la constante de Planck. La condition de
pr\'equantification s'\'ecrit  $p = 2 \pi \hbar n$, o\`u $n \in \mbox{\bf Z}$;
autrement dit,  $s = p/4\pi$ doit \^etre \'egal \`a  $\hbar n / 2$ (cf.
\cite{So}).  Alors, on obtient une {\em pr\'equantification par fusion}
\cite{So}.

Par ailleurs,
$(M,\Lambda) = (S^2,\sigma) \times (S^1,0)$ est int\'egr\'ee
par le groupo\"{\i}de symplectique
$(S^2,\sigma) \times (S^2,-\sigma) \times (T^\ast(S^1),-dx \wedge dy)$
qui est aussi pr\'equantifiable.

\noindent
3) $\cal F$ poss\`ede une composante de Reeb qui supporte un cycle
\'evanouissant non trivial et dans ce cas aucune structure de Poisson $\Lambda$
n'est int\'egrable d'apr\`es \cite{Hec}.}
\end{iexems}

\setcounter{ith}{0}
\setcounter{equation}{0}

\section{Vari\'et\'es de Poisson r\'eguli\`eres}

\subsection{Formes feuillet\'ees symplectiques}

\hspace{5mm}
Soit $(M,{\cal F})$ une vari\'{e}t\'{e} feuillet\'{e}e r\'eguli\`ere.
Soient
\( (\Omega^{\ast}(M),d) \)
le complexe de De Rham et
\( (\Omega^{\ast}(M,{\cal F}),d) \)
le sous-complexe des
{\it formes relatives} qui
s'annulent sur les feuilles de $\cal F$.
Le complexe quotient
\( (\Omega^{\ast}({\cal F}),d_{\cal F}) \)
est le complexe des
{\it formes feuillet\'{e}es} et sa cohomologie \( H^{\ast}({\cal F}) \)
est la {\it cohomologie feuillet\'{e}e de} $(M,{\cal F})$.
\vspace{1ex}

Soit
\( \nu^{\ast}{\cal F} \)
le fibr\'{e} conormal dont les sections sont les 1-formes relatives.
Le choix d'un suppl\'{e}mentaire $N{\cal F}$ de $T{\cal F}$ d\'efinit une
d\'{e}composition
\( T^{\ast}M = \nu^{\ast}{\cal F} \oplus T^{\ast}{\cal F} \)
qui induit des d\'{e}compositions
$\Omega^{r}(M) = \bigoplus_{p+q=r} \Omega^{p,q}(M)$ et
$d = d_{0,1} + d_{1,0} + d_{2,-1}$.
\vspace{1ex}

Si l'on consid\`ere $\Omega^{r}(M)$
comme un module filtr\'e de degr\'e filtrant $p$, on obtient
la {\em suite espectrale de Leray-Serre}
$ E_{2}^{p,q} \Rightarrow H^{p+q}(M) $ dont la diff\'erentielle
d'ordre $r$ sera not\'ee $d_r$.
La projection du complexe de De Rham sur le complexe des formes feuillet\'ees
se
restreint en un isomorphisme de complexes
$(\Omega^{0,\ast}(M),d_{0,1}) \cong (\Omega^\ast({\cal F}),d_{\cal F})$ et donc
le terme $E_1^{0,q} \cong H^q({\cal F})$.
\vspace{1ex}

Une 2-forme feuillet\'{e}e
\( \sigma \in \Omega^{2}({\cal F}) \)
est dite {\it symplectique} si
\( d_{\cal F} \sigma = 0 \) et
\( \bigwedge^{k}\sigma \)
est non nulle en tout point de $M$ o\`{u} $dim{\cal F} = 2k$.
D'apr\`es \cite{DH}, une telle forme feuillet\'{e}e d\'{e}finit une structure
de Poisson r\'eguli\`ere $\Lambda$ sur $M$ dont le feuilletage
caract\'{e}ristique est  $\cal F$.

\subsection{Invariants cohomologiques} \label{inv}

\hspace{5mm}
A la structure de Poisson r\'eguli\`ere
\( \Lambda = ({\cal F},\sigma) \),
on associe les \'{e}l\'{e}ments
\( [\sigma] = [\omega] \in E_{1}^{0,2} \)
et
\( d_{1}[\sigma] = [d_{1,0}\omega] \in E_{1}^{1,2} \)
de la suite spectrale de Leray-Serre,
o\`{u} $\omega$ est un repr\'{e}sentant pur de type $(0,2)$
de $\sigma$.
\vspace{1ex}

Si $\sigma$ poss\`ede un repr\'esentant pur $\omega$ tel que
$d_{1,0}\omega=0$, la forme volume des feuilles $v = \wedge^k \omega$ v\'erifie
aussi $d_{1,0}v = 0$ et il existe une m\'etrique riemannienne pour laquelle
les feuilles de $\cal F$ sont des sous-vari\'et\'es minimales.
On  dit que $\Lambda$ est {\em minimale} \cite{DH}. Le proc\'ed\'e de
{\em purification} de \cite{Su} (dont la preuve de la proposition~\ref{fmin}
sera un cas particulier) permet de montrer la proposition suivante:

\begin{prop}[\cite{DH}]  \label{min}
Une structure de Poisson r\'eguli\`ere $\Lambda$ est minimale si, et
seulement si, la classe $d_1[\sigma]$ est nulle. $\Box$
\end{prop}

Si $d_1[\sigma]=0$, on associe \`a $\Lambda$ des
nouvelles classes $[\sigma] \in E_2^{0,2}$ et
$d_2[\sigma] \in E_2^{2,1}$. Si $d_2[\sigma]=0$, on a des classes
$[\sigma] \in E_3^{0,2}$ et
$d_3[\sigma] \in E_3^{3,0}$.
En fait, d'apr\`es la proposition~\ref{min}, les classes $d_2[\sigma]$ et
$d_3[\sigma]$
appartiennent aux termes $E_1^{2,1}$ et
$E_2^{3,0}$ (cf. \S\ref{obst}).

\begin{prop}[\cite{DH},\cite{V1}]    \label{pres}
Une structure de Poisson r\'eguli\`ere $\Lambda$ est pr\'esymplectique
(i.e. $\sigma$ poss\`ede un repr\'esentant ferm\'e $\omega$) si, et seulement
si, les classes $d_1[\sigma]$, $d_2[\sigma]$ et $d_3[\sigma]$ sont nulles.
$\Box$
\end{prop}

\setcounter{th}{0}
\setcounter{equation}{0}

\section{Fibr\'es symplectiques} \label{fs}

\hspace{5mm}
Soit
\(\pi : M \rightarrow B \)
un fibr\'e symplectique. Dans ce paragraphe, on v\'erifiera que la structure
de Poisson associ\'ee $\Lambda = ({\cal F},\sigma)$ est minimale. Puis
on retrouvera le th\'eor\`eme suivant
\`a l'aide des autres obstructions $d_2[\sigma]$ et $d_3[\sigma]$:

\begin{th}[\cite{GLSW}]  \label{thGLSW}
Si la fibre $F$ et la base $B$ du fibr\'e symplectique sont simplement
connexes, alors $\sigma$ poss\`ede un repr\'esentant ferm\'e.
\end{th}

L'annulation de la classe $d_2[\sigma]$ est une cons\'equence de l'hypoth\`ese
sur la fibre $F$. Par ailleurs, l'ingr\'edient fondamental de la preuve
du th\'eor\`eme~\ref{thGLSW} consiste \`a montrer que le morphisme de
restriction
\( H^2(M) \rightarrow H^2(F) \)
est surjectif.
Pour cela, on prouve dans \cite{GLSW} que le cobord de la suite exacte
d'homotopie
\( \partial_\ast : \pi_3(B) \rightarrow \pi_2(F) \)
est d'ordre fini. On se propose d'utiliser ce fait pour
montrer que la classe $d_3[\sigma]$ est aussi nulle.

\subsection{Minimalit\'e des fibr\'es symplectiques} \label{Smin}

\hspace{5mm}
Soit $\omega$ un repr\'esentant pur de type $(0,2)$ de
$\sigma$ pour le choix d'un suppl\'ementaire de $T{\cal F}$. On va
d\'emontrer que la premi\`ere obstruction
\( d_1[\sigma] = [d_{1,0}\omega] \)
\`a l'existence d'un repr\'esentant ferm\'e est toujours nulle:

\begin{prop}  \label{fmin}
La structure de Poisson $\Lambda$ est minimale.
\end{prop}

\noindent
{\bf D\'emonstration}
Soient $\{U_i\}$ un recouvrement de $B$ form\'e d'ouverts qui trivialisent
$\pi$ et
\( \varphi_i : \pi^{-1}(U_i) \rightarrow U_i \times F \)
les cartes de trivialit\'e locale. Soit
\( p_2 : U_i \times F \rightarrow F \)
la projection. Si $\{\rho_i\}$ est une partition de l'unit\'e subordonn\'ee,
alors la 2-forme
\( \mu = \sum \pi^\ast(\rho_i) (p_2 \circ \varphi_i)^\ast (\sigma_F) \)
repr\'esente la 2-forme feuillet\'ee $\sigma$. Pour tout
couple $X_1$ et $X_2$ de champs tangents \`a $\cal F$, cette forme
v\'erifie

\begin{equation}    \label{noyau}
i_{\textstyle X_1} i_{\textstyle X_2} d\mu = 0
\end{equation}

\noindent
Puisque la 2-forme $\mu$ est non d\'eg\'en\'er\'ee sur le fibr\'e
vertical $T{\cal F}$, le
sous-fibr\'e horizontal
\( \{ \; Y \in TM \; / \; i_Y \mu \! \mid_{T{\cal F}} = 0 \; \} \)
est un suppl\'ementaire de $T{\cal F}$. Pour la d\'ecomposition
correspondante, la composante de type $(1,1)$ est nulle et donc
\( \mu = \mu^{0,2} + \mu^{2,0} \).
Si l'on pose $\omega = \mu^{0,2}$, ce nouveau repr\'esentant v\'erifie
$d_{1,0} \omega = 0$ d'apr\`es (\ref{noyau}). $\Box$

\subsection{Obstructions \`a l'existence d'un repr\'esentant ferm\'e}
\label{obst}

\hspace{5mm}
La 3-forme pure
\( d\omega = d_{2,-1} \omega \)
de type $(2,1)$ repr\'esente la deuxi\`eme obstruction
$d_2[\sigma] \in E_1^{2,1}$. Puisque la fibre $F$ est simplement connexe, ce
terme est nul d'apr\`es \cite{DH} (voir aussi le corollaire~\ref{corcoh}).
Donc $d_{2,-1}\omega$ poss\`ede une $d_{0,1}$-primitive $\eta$ de type $(2,0)$.

Si l'on d\'esigne encore par $\omega$ le nouveau repr\'esentant
$\omega - \eta$ de $\sigma$, alors la 3-forme
\( d\omega = d_{1,0} \eta \)
est {\em basique}, c'est-\`a-dire $d\omega$ se projette en une
3-forme ferm\'ee $\zeta$ sur $B$.

\begin{prop} Si la fibre $F$ et la base $B$ sont simplement connexes,
alors la troisi\`eme obstruction
\( [\zeta] \in H^3(B) \)
est nulle.
\end{prop}

\noindent
{\bf D\'emonstration}
On consid\`ere le diagramme commutatif suivant:

\begin{center}
\begin{picture}(220,49)(0,0)
\put(75,37){\vector(1,0){50}}
\put(75,5){\vector(1,0){50}}
\put(53,28){\vector(0,-1){14}}
\put(149,28){\vector(0,-1){14}}
\put(37,34){$\pi_3(B)$}
\put(133,34){$\pi_2(F)$}
\put(92,41){$\partial_\ast$}
\put(92,9){$d_3$}
\put(37,2){$H_3(B)$}
\put(133,2){$H_2(F)$}
\put(35,19){$H_3$}
\put(156,19){$H_2$}
\end{picture}
\end{center}

\noindent
Puisque $F$ est simplement connexe, le morphisme de Hurewicz $H_2$
est un isomorphisme et le morphisme $d_3$
est la {\em transgression} de la suite spectrale d'homologie
\cite{S}.
Si $B$ es simplement connexe, alors le
morphisme de Hurewicz $H_3$ est surjectif.
D'autre part, puisque le cobord
$\partial_\ast$ est d'ordre fini d'apr\`es \cite{GLSW},
il en est de m\^eme pour la transgression $d_3$. Donc le morphisme

\[ [z] \in H_3(B) \;
\stackrel{d_3}{\longmapsto} \;
[z_F] \in H_2(F) \;
\longmapsto \; \int_{\textstyle z_F} \sigma_F \in \mbox{\bf R} \]

\noindent
est nul. D'apr\`es la d\'efinition de transgression, le 2-cycle $z_F$ de $F$
est le bord d'une 3-cha\^{\i}ne $\tilde{z}$ de $M$ qui se projette sur le
3-cycle  $z$ de $B$. Alors, toute p\'eriode

\[ \int_{\textstyle z} \zeta =
\int_{\textstyle \tilde{z}} d\omega =
\int_{\textstyle z_F} \omega =
\int_{\textstyle z_F} \sigma_F \]

\noindent
est nulle; d'o\`u la proposition. $\Box$

\setcounter{th}{0}
\setcounter{equation}{0}

\section{Structure cohomologique des submersions} \label{coh}

\subsection{Th\'eor\`eme fondamental} \label{annul}

\hspace{5mm}
On commence par le rappel d'un
th\'eor\`eme de \cite{AH} sur la cohomologie d'une submersion surjective
$\pi : M \rightarrow B$:

\begin{th}                 \label{subm2}
Soit $\cal Q$ le faisceau image
r\'eciproque d'un faisceau ${\cal Q}_0$ de base $B$. Si les fibres $F_b$ sont
connexes et
$H^{q}(F_b;{\cal Q})=0$ en degr\'e $q < r$, alors
\vspace{1ex}

\noindent
i)
\( \;\; \pi^{\ast} : H^q(B;{\cal Q}_0) \rightarrow
H^{q}(M;{\cal Q}) \;\; \)  est isomorphisme en degr\'e $q < r$;
\vspace{1ex}

\noindent
ii) la suite
\( 0 \rightarrow
H^{r}(B;{\cal Q}_0) \stackrel{\pi^{\ast}}{\longrightarrow}
H^{r}(M;{\cal Q}) \stackrel{\rho}{\longrightarrow}
\prod_{b \in B} H^{r}(F_{b};{\cal Q}) \)
est exacte. $\Box$
\end{th}

Soit $\cal F$ le feuilletage d\'efini par $\pi$. Le
faisceau $\phi^p$ des germes de p-formes basiques est l'image r\'eciproque du
faisceau mou
\( \raisebox{-1.1ex}{$\stackrel{\textstyle \Omega}{\sim}$}^p(B) \)
des germes de p-formes sur $B$. D'apr\`es
\cite{V1}, on sait que
\( H^{q}(M;\phi^{p}) \cong  E_1^{p,q} \)
et l'on obtient le corollaire suivant:

\begin{cor}          \label{corcoh}
Dans les conditions du th\'eor\`eme \ref{subm2}, le terme
$E_1^{p,q} = 0$ pour tout $0 < q < r$ et le morphisme de restriction
$ \rho : E_1^{p,r} \longrightarrow \prod_{b \in B} H^{r}(F_{b};\phi^p) $
est injectif. $\Box$
\end{cor}

\subsection{Existence de repr\'esentants ferm\'es}    \label{repfer}

\hspace{5mm}
Le th\'eor\`eme~\ref{subm2} permet de d\'emontrer un crit\`ere d'existence de
repr\'esentants
ferm\'es qui g\'en\'eralise le th\'eor\`eme 1 de \cite{GLSW}:

\begin{th}   \label{subm1}
Soit
\( \pi : M \rightarrow B\)
une submersion surjective dont les fibres sont connexes et
cohomologiquement triviales en d\'egr\'e $q < r$.
Une $r$-forme feuillet\'ee ferm\'ee $\sigma$
poss\`ede un repr\'esentant ferm\'e si, et seulement si, il
existe une classe de cohomologie de De Rham
dont la restriction \`a chaque fibre $F_b$ est  \'egale \`a la
restriction de la classe
\( [\sigma] \in H^r({\cal F}) \).
\end{th}

\noindent
{\bf D\'emonstration}
On consid\`ere le diagramme de restriction:

\begin{center}
\begin{picture}(120,45)(0,0)
\put(42,32){\vector(1,-1){15}}
\put(42,37){\vector(1,0){50}}
\put(92,32){\vector(-1,-1){15}}
\put(0,34){$H^r(M)$}
\put(100,34){$H^r(\cal F)$}
\put(40,4){$\prod_{b \in B} H^r(F_b)$}
\put(39,20){$\rho$}
\put(89,20){$\rho_{\cal F}$}
\end{picture}
\end{center}

\noindent
Soit
\( [\mu] \in H^r(M) \)
une classe telle que
\( \rho([\mu]) = \rho_{\cal F} ([\sigma]) \).
Si l'on note $\overline{\mu}$ la classe feuillet\'ee de $\mu$, alors on a
\( \rho_{\cal F} ([ \sigma - \overline{\mu} ]) = 0 \).
Or, d'apr\`es le corollaire~\ref{corcoh},
le morphisme $\rho_{\cal F}$ est injectif
et donc la classe $[\sigma - \overline{\mu}] = 0$.
Soit $\chi$ un repr\'esentant
de la $d_{\cal F}$-primitive de $\sigma - \overline{\mu}$.
Alors, $\omega=\mu + d\chi$
est un repr\'esentant ferm\'e de $\sigma$; d'o\`u le th\'eor\`eme. $\Box$
\vspace{2ex}

L'exemple suivant montre que la trivialit\'e cohomologique des
fibres est essentielle dans la preuve du th\'eor\`eme~\ref{subm1}:

\begin{exem}
{\em Soit $\cal F$ le feuilletage de
\( \Delta = D^2 \times \lbrack 0, 1 \rbrack - \{(0,0)\} \)
d\'efini par l'\'equation $dt=0$. C'est le mod\`ele des {\em cycles
\'evanouissants
coh\'erents} \cite{Hec}. Le groupe
$$ H^2({\cal F}) =
{\cal C}^\infty (\rbrack 0, 1 \rbrack) / {\cal C}^\infty
(\lbrack 0, 1 \rbrack) $$
n'est pas nul, bien que les fibres sont
cohomologiquement triviales en degr\'e $2$.
Pour calculer la cohomologie feuillet\'ee, on d\'ecompose $\Delta$ en r\'eunion
de deux ouverts
\( \Delta_1 = D^2 \times \rbrack 0, 1 \rbrack \)
et
\( \Delta_2 = ( D^2 -\{0\} ) \times \lbrack 0, 1 \rbrack \)
dont l'intersection sera not\'ee $\Delta_0$. Soient ${\cal F}_1$, ${\cal F}_2$
et ${\cal F}_0$ les feuilletages induits. On
consid\`ere la {\em suite exacte de Mayer-Vietoris}

\[ 0 \rightarrow H^1({\cal F}) \rightarrow
H^1({\cal F}_1) \oplus H^1({\cal F}_2) \rightarrow H^1({\cal F}_0)
\rightarrow H^2({\cal F}) \rightarrow H^2({\cal F}_1) \oplus H^2({\cal F}_2)
\rightarrow \dots \]

\noindent
Puisque les feuilles de ${\cal F}_1$ sont contractiles,
le groupe \( H^q({\cal F}_1) = 0 \)
pour $q \geq 1$. Par ailleurs, $(\Delta_2,{\cal F}_2)$ se r\'etracte
par d\'eformation int\'egrable sur le bord
\( \partial \Delta_2 = S^1 \times \lbrack 0, 1 \rbrack \)
muni de la fibration en cercles.
Par int\'egration sur les fibres, il vient
\( H^1({\cal F}_2) \cong {\cal C}^\infty(\lbrack 0, 1 \rbrack) \)
et
\( H^2({\cal F}_2) = 0 \).
De m\^eme, on a
\( H^1({\cal F}_0) \cong {\cal C}^\infty(\rbrack 0, 1 \rbrack) \).
D'o\`u le calcul.}
\end{exem}

\setcounter{th}{0}
\setcounter{equation}{0}

\section{Submersions symplectiques} \label{Ssub}

\hspace{5mm}
Soit
\( \pi : M \rightarrow B \)
une submersion symplectique, i.e.
munie d'une forme feuillet\'ee symplectique $\sigma$
pour laquelle la structure de Poisson $\Lambda=({\cal F},\sigma)$ est minimale.

\subsection{P\'eriodes sph\'eriques}   \label{per}

\hspace{5mm}
Si $N$ d\'esigne le p\^ole nord, toute {\em sph\`ere tangente}
\( s : (S^{2},N) \rightarrow (F_b,x) \)
est contenue dans
un ouvert produit $U \times V$, o\`u $U$ est un voisinage de
$b$ dans $B$ et $V$ est un voisinage de $x$ dans la fibre $F_b$.
On obtient ainsi une
{\em d\'eformation transverse}
$$\begin{picture}(100,45)(0,0)
\put(40,27){\vector(1,-1){13}}
\put(40,32){\vector(1,0){40}}
\put(80,27){\vector(-1,-1){13}}
\put(0,29){$U \times S^2$}
\put(88,29){$M$}
\put(56,3){$U$}
\put(33,16){$p_1$}
\put(77,16){$\pi$}
\put(57,36){$D$}
\end{picture}$$
L'image r\'eciproque $D^\ast \sigma$ est une forme feuillet\'ee
ferm\'ee qui est repr\'esent\'ee par une 2-forme $\mu$ de type
$(0,2)$ pour la trivialisation canonique de $T(U \times S^2)$.
Par int\'egration sur les fibres de $p_1$,
on obtient une fonction diff\'erentiable
\( - \!\!\!\!\!  \int  \mu  \)
sur $U$. Les propri\'et\'es de $- \!\!\!\!\! \int$ impliquent que celle-ci
est ind\'ependante de la trivialisation de $T(U \times S^2)$. D'autre part,
\( - \!\!\!\!\!  \int  \mu  \)
ne d\'epend que des classes
d'homotopie des sph\`eres d'apr\`es le th\'eor\`eme de Stokes.
Ces fonctions d'aire d\'efinissent un
sous-groupo\"{\i}de $Per(\sigma)$ de
$B \times \mbox{\bf R}$ appel\'e le
{\em groupo\"{\i}de des p\'eriodes sph\'eriques de $\sigma$}.
Son image par la projection $p_2 : B \times \mbox{\bf R} \to \mbox{\bf R}$
engendre un sous-groupe $\cal P$ de {\bf R} appel\'e le {\em groupe des
p\'eriodes sph\'eriques de $\sigma$}.

La 3-forme $d_{1,0} \mu$ de type $(1,2)$ repr\'esente
l'image r\'eciproque de $d_1[\sigma]$. Les int\'egrales
\(  - \!\!\!\!\!  \int d_{1,0} \mu =
d( - \!\!\!\!\!  \int \mu) \)
d\'efinissent un sous-groupo\"{\i}de
\( Per_1(\sigma) \)
de $T^\ast B$ appel\'e le
{\em groupo\"{\i}de des p\'{e}riodes sph\'{e}riques d\'{e}riv\'{e}es
de $\sigma$}. Pour une construction g\'en\'erale, voir \cite{AH}.

\subsection{Submersions symplectiques}   \label{sub}

\hspace{5mm}
Le groupo\"{\i}de d\'eriv\'e d'une submersion symplectique est \'evidemment
nul. Cette condition va caract\'eriser certaines submersions symplectiques:

\begin{prop} Soient               \label{smin}
\( \pi : M \rightarrow B \)
une submersion surjective \`a fibres simplement connexes et
$\sigma$ une forme feuillet\'ee symplectique.
Si le groupo\"{\i}de d\'eriv\'e $Per_1(\sigma)$ est nul,
alors la submersion est symplectique. Si $H^3(B)$ est de plus nul,
alors $\sigma$ poss\`ede un repr\'esentant ferm\'e.
\end{prop}

\noindent
{\bf D\'emonstration}
Puisque le faisceau
\( \raisebox{-1.1ex}{$\stackrel{\textstyle \Omega}{\sim}$}^1(B) \)
est mou et les fibres $F_b$ sont simplement connexes, le morphisme de
restriction
\begin{center}
\( \rho : E_1^{1,2} = H^2(M;\phi^1) \longrightarrow \prod_{b \in B}
H^2(F_b;\phi^1) \)
\end{center}

\noindent
est injectif d'apr\`es le corollaire~\ref{corcoh}. Par ailleurs,
l'annulation de $Per_1(\sigma)$ signifie que les int\'egrales d'un
repr\'esentant de  $d_1[\sigma]$ sur les sph\`eres tangentes aux fibres sont
nulles. Donc les restrictions de $d_1[\sigma]$ aux fibres sont nulles.
Il s'ensuit que $d_1[\sigma]=0$ et donc la submersion est symplectique.
D'autre part, les conditions entra\^{\i}nent
l'annulation de $E_1^{2,1}$ et $E_2^{3,0} \cong H^3(B)$; d'o\`u la
deuxi\`eme affirmation. $\Box$

\subsection{La condition d'int\'egrabilit\'e: d\'emonstration du
th\'eor\`eme 1}   \label{th1}

\hspace{5mm}
Soit
\( \pi : M \rightarrow B \)
une submersion symplectique \`a fibres connexes et simplement connexes telle
que
$H^3(B;\mbox{\bf Z})=0$. Pour d\'emontrer le th\'eor\`eme~\ref{condint},
on supposera que $\cal P$ est un sous-groupe discret de {\bf R} et
l'on proc\'edera en deux \'etapes:
\vspace{1ex}

\noindent
1) {\bf Int\'{e}gration cohomologique}: on montrera que la classe "r\'eelle"
$[\sigma] \in H^2({\cal F})$
s'in\-t\`egre
en une classe "enti\`ere"
$ \nu \in H^2(M;{\cal P}) $. Le th\'eor\`eme~\ref{subm2} ram\`enera
l'int\'egration globale
\`a l'int\'egration fibre \`a fibre.
\vspace{2ex}

\noindent
2) {\bf Int\'{e}gration diff\'erentiable}: on v\'erifiera que
la classe $\nu$ est repr\'esent\'ee par un cocycle sur
$M$ \`a valeurs dans $\mbox{\bf R}/{\cal P}$ qui d\'efinira la
pr\'equantification.
\vspace{2ex}

\noindent
1) Le faisceau $\phi^0$ des germes de fonctions basiques est
l'image r\'ecripoque du faisceau
\( \raisebox{-1.1ex}{$\stackrel{\textstyle {\cal C}}{\sim}$}^\infty(B) \)
des germes de fonctions sur $B$. Soient $\cal Q$ et ${\cal Q}_0$ les
quotients des faisceaux $\phi^0$ et
\( \raisebox{-1.1ex}{$\stackrel{\textstyle {\cal C}}{\sim}$}^\infty(B)  \)
par les faisceaux constants de fibre $\cal P$ (que l'on notera de la m\^eme
fa\c{c}on). Les suites exactes correspondantes induisent des suites exactes
longues de cohomologie

\begin{equation}    \label{suite1}
\raisebox{3ex}{
\begin{array}[t]{c}
H^2(B;{\cal P}) \\
\; \; \; \downarrow  \pi^\ast \\
H^2(M;{\cal P})
\end{array}
\begin{array}[t]{c}
\stackrel{i^\ast_0}{\longrightarrow} \\  \\
\stackrel{i^\ast}{\longrightarrow}
\end{array}
\begin{array}[t]{c}
H^2(B;\raisebox{-1.1ex}{$\stackrel{\textstyle {\cal C}}{\sim}$}^\infty(B)) \\
\; \; \; \downarrow \pi^\ast \\
H^2(M;\phi^0)
\end{array}
\begin{array}[t]{c}
\stackrel{j^\ast_0}{\longrightarrow} \\ \\
\stackrel{j^\ast}{\longrightarrow}
\end{array}
\begin{array}[t]{c}
H^2(B;{\cal Q}_0) \\
\; \; \; \downarrow \pi^\ast \\
H^2(M;{\cal Q})
\end{array}
\begin{array}[t]{c}
\stackrel{\delta_0}{\longrightarrow} \\ \\
\stackrel{\delta}{\longrightarrow}
\end{array}
\begin{array}[t]{c}
H^3(B;{\cal P}) \\
\; \; \; \downarrow \pi^\ast \\
H^3(M;{\cal P})
\end{array}}
\end{equation}

\noindent
On remarque que:
\vspace{1ex}

\noindent
i) le groupe $H^2(M;\phi^0)$ est isomorphe au groupe
$H^2({\cal F})$ d'apr\`es \cite{V1};
\vspace{1ex}

\noindent
ii) le cobord $\delta_0$ est un isomorphisme, car le faisceau
\( \raisebox{-1.1ex}{$\stackrel{\textstyle {\cal C}}{\sim}$}^\infty(B) \)
est mou;
\vspace{1ex}

\noindent
iii) le groupe $H^2(B;{\cal Q}_0)=0$, car $H^3(B;{\cal P})=0$.
\vspace{2ex}

Les inclusions des fibres $F_b$
dans $M$ induisent des morphismes de restriction

\begin{equation}     \label{rest}
\raisebox{-5ex}{
\begin{picture}(305,52)(0,0)
\put(80,44){\vector(1,0){40}}
\put(80,6){\vector(1,0){40}}
\put(12,41){$H^{2}(M;{\cal P})$}
\put(40,33){\vector(0,-1){15}}
\put(2,3){$\prod_{b \in B} H^{2}(F_b;{\cal P})$}
\put(208,44){\vector(1,0){40}}
\put(208,6){\vector(1,0){40}}
\put(138,41){$H^{2}(M;\phi^{0})$}
\put(265,41){$H^{2}(M;{\cal Q})$}
\put(166,33){\vector(0,-1){15}}
\put(291,33){\vector(0,-1){15}}
\put(127,3){$\prod_{b \in B} H^{2}(F_b;\phi^{0})$}
\put(252,3){$\prod_{b \in B} H^{2}(F_b;{\cal Q})$}
\put(90,11){$\{i^{\ast}_b\}$}
\put(96,48){$i^{\ast}$}
\put(224,48){$j^{\ast}$}
\put(218,11){$\{j^{\ast}_b\}$}
\put(48,23){$\rho$}
\put(174,23){$\rho$}
\put(299,23){$\rho$}
\end{picture}}
\end{equation}
\vspace{0.5ex}

\noindent
et l'on a une factorisation

\begin{equation}    \label{fact}
\raisebox{-5ex}{
\begin{picture}(135,50)(0,0)
\put(62,33){\vector(1,-1){17}}
\put(62,41){\vector(1,0){50}}
\put(92,19){\vector(1,1){15}}
\put(6,37){$H^{2}(F_b;{\cal P})$}
\put(120,37){$H^{2}(F_b;\phi^{0})$}
\put(70,4){$H^{2}(F_b)$}
\put(54,21){$h^\ast_b$}
\put(108,21){$k^\ast_b$}
\put(82,45){$i^\ast_b$}
\end{picture}}
\end{equation}
\vspace{0.5ex}

Par int\'{e}gration de $\sigma$ sur les sph\`{e}res tangentes aux fibres
(qui sont simplement connexes),
on obtient des classes "enti\`eres"
$$ \nu_{b} \in Hom(\pi_2(F_b),{\cal P}) = H^{2}(F_u;{\cal P}) $$
qui int\`egrent les classes "r\'eelles" $[\sigma_b] \in H^2(F_u)$,
c'est-\`a-dire
\( h^\ast_b (\nu_b) = [\sigma_b]\). On en d\'eduit que
$ \rho([\sigma]) = \{ i^{\ast}_b (\nu_{b}) \} $ d'apr\`es
la factorisation~(\ref{fact}) et donc
$$ \rho(j^{\ast}([\sigma])) = \{ j^\ast_b \} (\rho(\sigma)) =
\{ j^{\ast}_b {\scriptstyle \circ} i^{\ast}_b (\nu_b) \} = 0 $$
car le diagramme~(\ref{rest}) est commutatif.
Il s'ensuit que la classe $j^\ast[\sigma]$ appartient au groupe
$H^2(B;{\cal Q}_0)$ qui est le noyau du morphisme $\rho$
d'apr\`es le th\'eor\`eme~\ref{subm2}. Or ce groupe est nul et donc la classe
$j^\ast[\sigma]= 0$. L'exactitude de la suite~(\ref{suite1}) implique que la
classe $[\sigma] \in H^2({\cal F})$ se remonte en une classe
$\nu \in H^2(M;{\cal P})$.
\vspace{1ex}

\noindent
2) A la suite exacte de groupes
$ 0 \rightarrow {\cal P}
\rightarrow \mbox{\bf R}
\rightarrow \mbox{\bf R}/{\cal P}
\rightarrow 0 $,
on associe la suite exacte de faisceaux de germes de fonctions
$$ 0 \longrightarrow
{\cal P} =
\raisebox{-1.1ex}{$\stackrel{\textstyle {\cal C}}{\sim}$}^\infty(M,{\cal P})
\longrightarrow
\raisebox{-1.1ex}{$\stackrel{\textstyle {\cal C}}{\sim}$}^\infty(M)
\longrightarrow
\raisebox{-1.1ex}{$\stackrel{\textstyle {\cal C}}{\sim}$}^\infty(M,\mbox{\bf R}
/{\cal P})
\longrightarrow 0 $$
Le faisceau
$\raisebox{-1.1ex}{$\stackrel{\textstyle {\cal C}}{\sim}$}^\infty(M)$
\'etant mou, le cobord de la suite exacte longue de cohomologie
$$ \delta :
H^1(M;
\raisebox{-1.1ex}{$\stackrel{\textstyle {\cal C}}{\sim}$}^\infty(M,\mbox{\bf
R}/
{\cal P}))  \longrightarrow  H^2(M;{\cal P})
$$
est un isomorphisme. La classe $\tau = \delta^{-1}\nu$ est repr\'esent\'ee
par  un cocycle sur $M$ \`a valeurs dans $\mbox{\bf R} / {\cal P} \cong
\mbox{\bf R} / p \mbox{\bf Z}$.  Ce cocycle d\'efinit un $S^1$-fibr\'e
principal $q: E \to M$.  La longeur des fibres est \'egale \`a $p$ (except\'e
le cas  $p=0$ o\`u le fibr\'e est trivial).
\vspace{1ex}

Le morphisme naturel
\( h^\ast : H^2(M;{\cal P}) \rightarrow H^2(M) \)
envoie $\nu$ sur la {\em classe de Chern}, i.e. la classe de la courbure
$\omega_0$ d'une connexion $\theta_0$ sur $E$. Puisque
$i^\ast\nu = [\sigma]$, la classe feuillet\'ee $\overline{\omega}_0$ de
$\omega_0$ est cohomologue \`a $\sigma$. Si $\chi$ est un repr\'esentant
d'une $d_{\cal F}$-primitive de $\sigma - \overline{\omega}_0$,
la courbure $\omega = \omega_0 + d\chi$ de la connexion
$\theta = \theta_0 + q^\ast \chi$ repr\'esente $\sigma$.
\vspace{2ex}

R\'eciproquement, soit
\( q : E \rightarrow M \)
une pr\'equantification munie d'une
connexion de courbure $\omega$. Le groupe des p\'eriodes $Per(\omega)$
est alors un sous-groupe discret de {\bf R}.
Puisque $\omega$ repr\'esente $\sigma$, $\cal P$ est un sous-groupe de
$Per(\omega)$; d'o\`u le th\'eor\`eme~\ref{condint}. $\Box$

\begin{rem}     \label{disden}
{\em Dans les conditions du th\'eor\`eme~\ref{condint}, la forme
feuillet\'ee symplectique $\sigma$ poss\`ede toujours des repr\'esentants
ferm\'es d'apr\`es la proposition~\ref{smin}. Mais l'int\'egration
cohomologique
fournit un bon repr\'esentant ferm\'e dont le groupe des
p\'eriodes est \'egal \`a $\cal P$.
L'exemple suivant montre que cela n'est pas vrai en g\'en\'eral. Soit
$ \pi : M = \mbox{\bf R}^2 \times S^2 \times S^2 \to B = S^2 \times S^2$
la fibration triviale. Soient $p_1$ et $p_2$ les projections de
$B = S^2 \times S^2$ sur chacun des facteurs et $\rho$ un nombre irrationnel.
Si l'on note $v$ la forme volume canonique sur $S^2$, la 2-forme ferm\'ee
$\omega = dp \wedge dq  + \pi^\ast(p_1^\ast v + \rho p_2^\ast v)$
repr\'esente une forme feuillet\'ee pour laquelle
$\pi$ devient une fibration symplectique. Le groupe $\cal P$ est nul, mais les
p\'eriodes sph\'eriques de $\omega$ sont partout denses dans {\bf R}.}
\end{rem}

\subsection{Exemples}      \label{exe2}

\noindent
1) Soit
\( \pi : M \rightarrow B \)
une submersion \`a fibres connexes et simplement connexes munie d'une forme
feuillet\'ee symplectique $\sigma$. Si les fibres sont en outre asph\'eriques,
la classe $[\sigma]$ est  nulle d'apr\`es le corollaire~\ref{corcoh}. Il
s'ensuit que la submersion est symplectique.  De plus, celle-ci est
pr\'equantifi\'ee par un fibr\'e trivial en cercles.  \vspace{1ex}

\noindent
2) Soit $N$ le p\^ole nord de la sph\`ere $S^2$.
Soit $\cal F$ le feuilletage de
\( M = S^{2} \times {\bf R} - \{ N \} \times [0,1] \)
d\'efini par l'\'equation $dt=0$.
En multipliant la forme volume normalis\'ee de $S^2$ par une
fonction positive convenable, on peut obtenir une forme volume des
feuilles $\omega$ dont la fonction d'aire
$ - \!\!\!\!\! \int_{S^2} \omega \in {\cal C}^\infty(\mbox{\bf R} - [0,1]) $
est \'egale \`a $1$ si $t < 0$ et \`a $2$ si $t > 1$.
La forme feuillet\'ee correspondante $\sigma$ d\'efinit
une structure de submersion symplectique sur $M$ telle que ${\cal P} =
\mbox{\bf Z}$.
La suite spectrale de \v{C}ech (cf. \cite{G}) attach\'ee au recouvrement
$ M_1 = S^2 \times \rbrack -\infty , 0 \lbrack $,
$ M_2 = S^2 \times \rbrack 1 , + \infty \lbrack $ et
$ M_3 = ( S^2 - \{ N \} ) \times \mbox{\bf R} $
permet de calculer la cohomologie feuillet\'ee. En effet,
le groupe $H^2({\cal F})$ se r\'eduit au terme $E_2^{0,2}$ de cette suite
spectrale.
Si l'on note ${\cal F}_i$ le feuilletage induit sur $M_i$, il s'ensuit:

\[ H^2({\cal F}) =
H^2({\cal F}_1) \oplus H^2({\cal F}_2) \oplus H^2({\cal F}_3)  =
H^2({\cal F}_1) \oplus H^2({\cal F}_2) \]

\noindent
o\`u $ H^2({\cal F}_3) = 0 $
car les feuilles sont contractiles. Par int\'egration sur
les fibres, il vient
$$ H^2({\cal F}) = {\cal C}^\infty ( \rbrack -\infty , 0 \lbrack ) \oplus
{\cal C}^\infty ( \rbrack 1 , + \infty \lbrack ) =
{\cal C}^\infty(\mbox{\bf R}-[0,1]) $$
De mani\`ere analogue, on montre que $H^{2}(M;\mbox{\bf Z}) = \mbox{\bf Z}
\oplus \mbox{\bf Z}$.
L'image de la classe $\nu=(1,2)$ par le morphisme
\( i^\ast : H^{2}(M;\mbox{\bf Z}) \to H^{2}({\cal F}) \)
est la classe $[\sigma]$ identifi\'ee \`a la fonction d'aire.
En restriction \`a chaque fibre, la pr\'equantification globale
induit la fibration triviale si $0 \leq t \leq 1$,
les fibrations de Hopf
\( S^3 \rightarrow S^2 \)
si $t < 0$ et
\( \mbox{\bf R}P^3 \rightarrow S^2 \)
si $t > 1$.
\vspace{1ex}

\noindent
3) Soit $\rho$ un nombre irrationnel. On peut remplacer la fonction
choisie dans l'exemple 2 de fa\c{c}on \`a obtenir une fonction d'aire
\'egale \`a $1$ pour $t<0$ et \`a $\rho$ pour $t>1$.
Dans ce cas, $\cal P$ est un sous-groupe partout dense de {\bf R}.

\subsection{Quelques remarques sur la pr\'equantification} \label{preq}

i) Si l'on cherche \`a pr\'equantifier les submersions
symplectiques par des fibr\'es en droites complexes,
la condition d'int\'egrabilit\'e est plus restrictive:
le groupe $\cal P$ doit \^etre un sous-groupe de {\bf Z} (pour le choix de
la  longeur $1$). Alors, le
th\'eor\`eme~\ref{condint} fournit un
repr\'esentant ferm\'e $\omega$ de $\sigma$ dont les
p\'eriodes sont enti\`eres. C'est le cas dans l'exemple 2 de \S\ref{exe2}.
Si
\(\; ({\cal L}^\ast(M,\Lambda),\partial) \; \)
est le {\em complexe de Lichn\'erowicz-Poisson} (cf. \cite{V2}) et
\( \# : (\Omega^\ast(M),d) \rightarrow  ({\cal L}^\ast(M,\Lambda),\partial) \)
est le morphisme de complexes induit par
$\Lambda$, on a la relation $\omega^\# = \Lambda$.
C'est un cas particulier de la condition de quantification de \cite{V2}.

La notion de pr\'equantification utilis\'ee dans ce travail est analogue \`a
celle de \cite{TW}. Dans ce travail, on trouvera
une bonne discussion sur les diff\'erentes notions de
pr\'equantification dans le cas classique.
\vspace{1ex}

\noindent
iii) Si $Per(\sigma)$ est un sous-groupo\"{\i}de
de Lie plong\'{e} de $B \times \mbox{\bf R}$ \'etal\'e sur $B$, alors
le quotient
\( {\cal G} = B \times \mbox{\bf R} / Per(\sigma) \)
est un groupo\"{\i}de de Lie. En proc\'edant comme dans
\cite{AH}, on construit un
fibr\'e principal
de groupo\"{\i}de structural $\cal G$ muni d'une connexion partielle
de courbure $\sigma$. L'exemple 3 de \S\ref{exe2} est "pr\'equantifiable" dans
ce sens large.

\setcounter{th}{0}
\setcounter{equation}{0}

\section{Pr\'equantification au sens de Weinstein}
\label{SpreW}

\hspace{5mm}
Le but de ce paragraphe est de montrer que l'emploi de l'int\'egration
symplectique est naturel si l'on cherche \`a pr\'equantifier les vari\'et\'es
de Poisson. Tout d'abord, le probl\`eme de la pr\'equantification deviendra le
probl\`eme de la construction  d'extensions centrales de groupo\"{\i}des
symplectiques (cf. \cite{WX}). On explicitera l'analogie  avec les r\'esultats
de \cite{TW} sur les extensions centrales de groupes. D'autre part,  {\em
l'int\'egration de Poisson} \cite{DH} developpe le feuilletage
caract\'eristique
en  une submersion surjective \`a fibres connexes et simplement connexes. Le
passage \`a cette  int\'egration interm\'ediaire ram\`enera l'\'etude
cohomologique des vari\'et\'es de Poisson  r\'eguli\`eres \`a celle des
submersions \`a fibres connexes simplement connexes.  En particulier, les
obstructions cohomologiques dispara\^{\i}tront  pour les vari\'et\'es de
Poisson
minimales.

\subsection{Extensions de groupo\"{\i}des symplectiques:
d\'emonstration du th\'eor\`eme~2}   \label{extg}

\hspace{5mm}
Soit $(M,\Lambda)$ une vari\'et\'e de Poisson pr\'equantifiable au sens
de Weinstein.
On se propose de d\'emontrer la premi\`ere partie du th\'eor\`eme~\ref{canon},
\`a savoir qu'il
existe une seule pr\'equantification
\( q : (E,\theta) \rightarrow (\Gamma,\eta) \)
triviale en restriction \`a l'espace des unit\'es $M$.
Cela signifie que l'holonomie du fibr\'e plat induit sur $M$ est triviale.
Les fibres de la projection source $\alpha$ sont connexes et
simplement connexes.
D'apr\`es le  th\'eor\`eme~\ref{subm2}, le morphisme
$ \alpha^\ast : H^1(M;S^1) \rightarrow H^1(\Gamma;S^1) $
est un isomorphisme dont l'inverse est
induit par l'inclusion $ \varepsilon : M \to \Gamma $.
Les pr\'equantifications de
$(\Gamma,\eta)$ sont classifi\'ees par
$H^1(\Gamma;S^1) \cong Hom(\pi_1(\Gamma),S^1)$ (voir \cite{TW}) et donc par
$H^1(M;S^1) \cong Hom(\pi_1(M),S^1)$. Une classe dans ce groupe est la
diff\'erence des morphismes d'holonomie des fibr\'es plats induits par
deux pr\'equantifications. On en d\'eduit qu'une
pr\'equantification de $(\Gamma,\eta)$ triviale en restriction \`a $M$
est unique \`a \'equivalence pr\`es.
\vspace{1ex}

D'apr\`es \cite{WX}, l'espace total $E$ est muni d'une structure canonique
de groupo\"{\i}de de Lie qui en fait une extension de $\Gamma$ pas $S^1$. Pour
retrouver cette structure, on compl\`ete  l'int\'egration symplectique en un
diagramme

\begin{equation}
\raisebox{-4ex}{
\begin{picture}(220,50)(0,0)
\put(53,33){\vector(2,-1){35}}
\put(50,41){\vector(1,0){50}}
\put(120,41){\vector(1,0){50}}
\put(167,33){\vector(-2,-1){35}}
\put(175,33){\vector(-2,-1){35}}
\put(107,33){\vector(0,-1){20}}
\put(113,33){\vector(0,-1){20}}
\put(6,37){$M \times S^1$}
\put(180,37){$\Gamma$}       \label{ext1}
\put(107,37){$E$}
\put(105,0){$M$}
\put(55,18){$p_{1}$}
\put(97,24){$\widehat{\alpha}$}
\put(118,24){$\widehat{\beta}$}
\put(141,28){$\alpha$}
\put(162,16){$\beta$}
\put(142,46){$q$}
\put(73,46){$i$}
\end{picture}}
\end{equation}
\vspace{0.5ex}

Soit ${\cal H}_\Gamma$ l'alg\`ebre de Lie des champs hamiltoniens
$X_{f {\scriptstyle \circ} \alpha}$ d\'efinis par
$$i_{\textstyle X_{f {\scriptstyle \circ} \alpha}} \eta = - \alpha^\ast df$$
o\`u $f \in {\cal C}^\infty(M)$. C'est une sous-alg\`ebre de Lie de
l'alg\`ebre de Lie ${\cal L}_\Gamma$ des {\em champs invariants \`a gauche} $X$
(cf. \cite{AD}) d\'efinis par $i_X \eta = - \alpha^\ast \mu$, o\`u
$\mu \in \Omega^1(M)$.

Soient $Z$ le champ
fondamental de l'action de $S^1$ sur $E$ et
$\widetilde{X}_{f {\scriptstyle \circ} \alpha}$
le rel\`evement horizontal du champ hamiltonien $X_{f {\scriptstyle \circ}
\alpha}$.
Les champs $Y_{f {\scriptstyle \circ} \alpha}$ d\'efinis par
$$ Y_{\textstyle f {\scriptstyle \circ} \alpha} =
\widetilde{X}_{\textstyle f {\scriptstyle \circ} \alpha}
+ q^\ast(f {\scriptstyle \circ} \alpha) Z$$
v\'erifient:
\vspace{1ex}

\noindent
i) $Y_{f {\scriptstyle \circ} \alpha}$ est un champ invariant par l'action de
$S^1$ qui se projette sur le champ $X_{f {\scriptstyle \circ} \alpha}$;
\vspace{1ex}

\noindent
ii) $L_{\textstyle Y_{f {\scriptstyle \circ} \alpha}} \theta = 0$;
\vspace{1ex}

\noindent
iii) $[Y_{{f_1 {\scriptstyle \circ} \alpha}},
Y_{{f_2 {\scriptstyle \circ} \alpha}}] =
Y_{\{f_1,f_2\} {\scriptstyle \circ} \alpha}$ pour tout couple de
fonctions  $f_1,f_2 \in {\cal C}^\infty(M)$.
\vspace{1ex}

\noindent
Ces champs forment une sous-alg\`ebre de Lie ${\cal H}_E$ de $\mbox{\bf X}(E)$.
De fa\c{c}on pr\'ecise, on obtient une extension centrale
$$ 0 \longrightarrow
\mbox{\bf R}
\longrightarrow
{\cal H}_E
\stackrel{q_\ast}{\longrightarrow}  {\cal H}_\Gamma
\longrightarrow 0$$

A l'alg\`ebre de Lie ${\cal H}_\Gamma$ des champs hamiltoniens
$Y_{f {\scriptstyle \circ} \alpha}$, on associe le faisceau
$\widetilde{\cal H}_\Gamma$ des germes correspondants. C'est un
{\em faisceau de d\'efinition} \cite{AD} du groupo\"{\i}de symplectique
$(\Gamma,\eta)$. Le faisceau d'alg\`ebre de Lie
$\widetilde{\cal H}_E$ associ\'e \`a l'extension v\'erifie:
\vspace{1ex}

\noindent
1) $\widetilde{\cal H}_E$ {\em s\'epare les fibres de $\widehat{\alpha}$} au
sens de \cite{AD}, car le faisceau d'alg\`ebre de Lie
$\widetilde{\cal H}_\Gamma$ s\'epare celles de $\alpha$;
\vspace{1ex}

\noindent
2) les $\widehat{\beta}$-fibres sont les
orbites transitives de l'action de $\widetilde{\cal H}_E$ car les
$\beta$-fibres  sont celles de l'action de $\widetilde{\cal H}_\Gamma$;
\vspace{1ex}

\noindent
3) le flot d'un champ
$Y_{f \circ \alpha}$ est form\'e de diff\'eomorphismes locaux d\'efinis sur des
ouverts satur\'es pour $\widehat{\alpha}$.
\vspace{1ex}

\noindent
D'apr\`es \cite{AD}, les flots des champs
$Y_{f \circ \alpha}$ forment un pseudo-groupe qui d\'efinit une structure de
groupo\"{\i}de de Lie sur $E$. La suite
exacte~(\ref{ext1}) devient alors une extension de $\Gamma$ par $S^1$.

\begin{rem}
{\em On dira qu'un champ $X$ est une {\em symm\'etrie infinit\'esimale} de
$(\Gamma,\eta)$ si $L_X \eta = 0$. Un champ invariant \`a gauche est une
symm\'etrie  infinit\'esimale si, et seulement si, il est localement
hamiltonien. Puisque l'alg\`ebre de  symm\'etries infinit\'esimales ${\cal
H}_\Gamma$ engendre un faisceau de d\'efinition de $\Gamma$,  on dira que
$\Gamma$ est un {\em groupo\"{\i}de de symm\'etries de $(\Gamma,\eta)$}.
D'autre part, soit $J$ l'application qui, \`a tout champ  $X_{f {\scriptstyle
\circ} \alpha} \in {\cal H}_\Gamma$,  associe la fonction $f {\scriptstyle
\circ} \alpha \in {\cal C}^\infty(\Gamma)$.  C'est une {\em application moment}
pour l'action \`a droite de $\Gamma$ sur $(\Gamma,\eta)$.  Cela explicite
l'analogie avec la construction   d'extensions centrales de groupes d'apr\`es
\cite{TW}. La description dans  \cite{WX} du cocycle de
l'extension~(\ref{ext1})
pr\'ecise cette analogie.}
\end{rem}

\subsection{Int\'egration de Poisson}

\hspace{5mm}
Soit $\Lambda_0=({\cal F}_0,\sigma_0)$ une structure de Poisson
r\'eguli\`ere sur une vari\'et\'e $M_0$ (o\`u l'on modifie les notations
ci-dessus par l'adjonction d'un indice $0$).

Le {\em groupo\"{\i}de d'homotopie} $\Pi_{1}({\cal F}_{0})$ de ${\cal F}_0$ est
le quotient de l'espace des chemins contenus dans les feuilles de
${\cal F}_{0}$ (muni de la topologie compact-ouvert $C^{\infty}$) par la
relation d'homotopie dans les feuilles de ${\cal F}_{0}$. C'est un
groupo\"{\i}de de Lie dont l'espace total $M$ est
s\'{e}par\'{e} si tout cycle \'{e}vanouissant de
${\cal F}_{0}$ est trivial (cf. \cite{DH}).
Les projections source $\alpha_{0}$ et but $\beta_{0}$ d\'{e}finissent
un m\^{e}me feuilletage image r\'{e}ciproque $\cal F$ de ${\cal F}_{0}$
dont la trace sur $M_0$ est \'egale \`a ${\cal F}_0$.

La structure de Poisson $\Lambda_{0}$ se rel\`{e}ve en une structure de
Poisson $\Lambda$ sur $M$ qui en fait un {\em groupo\"{\i}de de Poisson}
au sens de \cite{W} (voir \cite{DH}). Cette structure est d\'{e}termin\'{e}e
par le feuilletage $\cal F$ et la forme feuillet\'{e}
\( \sigma = \alpha_{0}^{\ast}\sigma_{0} - \beta_{0}^{\ast}\sigma_{0} \).

\begin{prop}[\cite{DH}]
Si la vari\'et\'e de Poisson $(M_0,\Lambda_0)$ est minimale (resp.
totalement asph\'erique), alors l'int\'egration de Poisson $(M,\Lambda)$
est pr\'esymplectique (resp. exacte) et donc $(M_0,\Lambda_0)$ est
int\'egrable.
\end{prop}

\noindent
{\bf D\'emonstration}
Les classes
$d_1[\sigma]$, $d_2[\sigma]$ et $d_3[\sigma]$
appartiennent aux termes $E_1^{1,2}$, $E_2^{2,1}$ et $E_3^{3,0}$ de la suite
spectrale de Leray-Serre de la paire $(M,M_0)$, car la forme feuillet\'ee
symplectique
\( \sigma = \alpha_{0}^{\ast}\sigma_{0} - \beta_{0}^{\ast}\sigma_{0} \)
s'annule en restriction \`a $M_0$.

Si $\Lambda_0$ est minimale, la classe $d_1[\sigma_0]$ est nulle et donc il en
est de m\^eme pour la classe $d_1[\sigma]$.
Puisque les fibres de $\alpha_0$ sont connexes
et simplement connexes, la version relative du corollaire~\ref{corcoh}
implique que les classes $d_2[\sigma]$ et $d_3[\sigma]$ sont aussi nulles.
Bref, $\sigma$ poss\`ede un repr\'esentant ferm\'e $\omega$.

Si $\Lambda_0$ est totalement asph\'erique, les feuilles de
${\cal F}_0$ (et donc les fibres de $\alpha_0$) sont asph\'eriques.
Dans ce cas, $E_1^{0,2}=0$ d'apr\`es la version relative du
corollaire~\ref{corcoh}. En
par\-ticulier, la classe $[\sigma]$ est nulle et donc $\sigma$ poss\`ede un
repr\'esentant exact $\omega$.

D'autre part, soit $L$ la forme de Liouville du
fibr\'e conormal $p : \Gamma = \nu^\ast{\cal F} \rightarrow M$. La 2-forme
\( \eta = p^\ast \omega - dL \)
est symplectique et le groupo\"{\i}de symplectique $(\Gamma,\eta)$
r\'ealise l'int\'egration symplectique de $(M_0,\Lambda_0)$. $\Box$
\vspace{2ex}

Par ailleurs, l'\'ecriture de $\sigma$ permet de d\'ecrire leurs
p\'eriodes sph\'eriques en termes de celles de $\sigma_0$:

\begin{lem}[\cite{AH}]   \label{comp}
Le groupo\"{\i}de $Per(\sigma)$ est l'image r\'eciproque par
$\alpha_0$ et $\beta_0$ du groupo\"{\i}de $Per(\sigma_0)$ et donc
${\cal P} = {\cal P}_0$. $\Box$
\end{lem}

En proc\'edant comme dans la preuve du th\'eor\`eme~\ref{condint}
(voir aussi \cite{AH}), on d\'emontre la condition de pr\'equantification
suivante:

\begin{th}    \label{Poisson}
Soit $(M_0,\Lambda_0)$ une vari\'et\'e de Poisson r\'eguli\`ere
dont tout cycle \'evanouis\-sant est trivial. L'int\'egration de Poisson
$(M,\Lambda)$ est pr\'equantifiable si, et seulement si,
${\cal P}_0$ est un sous-groupe discret de {\bf R}. $\Box$
\end{th}

\begin{rem}
{\em (1) La condition sur les cycles \'evanouissants de ${\cal F}_0$ implique
que l'int\'egration de
Poisson est s\'epar\'ee ce qui permet d'int\`egrer la classe $[\sigma]$.
\vspace{1ex}

\noindent
(2) La pr\'equatification de $(M,\Lambda)$
est obtenue \`a l'aide d'un cocycle dont la classe $\nu$
int\`egre la classe $[\sigma]$. Puisque cette classe est relative, la classe
$\nu$ est aussi relative et donc le $S^1$-fibr\'e principal induit au-dessus
de $M_0$ est trivial. D'autre part, les pr\'equantifications triviales en
restriction \`a $M_0$ sont  classifi\'ees par le groupe de cohomologie relative
$H^1(M,M_0;S^1)$. Puisque les fibres de  $\alpha_0$ sont connexes et simplement
connexes, ce groupe est nul d'apr\`es la version  relative du
th\'eor\`eme~\ref{subm2}. Bref, le th\'eor\`eme~\ref{Poisson} fournit la seule
pr\'equantification de $(M,\Lambda)$ qui est triviale en restriction \`a $M_0$.
\vspace{1ex}

\noindent
(3) L'espace total $E$ ne peut pas \^etre muni d'une structure de
groupo\"{\i}de, car les champs hamiltoniens $X_{f {\scriptstyle \circ}
\alpha_0}$ ne sont pas des symm\'etries  infinit\'esimales de $(M,\Lambda)$.
N\'eanmoins, on retrouve la situation du   th\'eor\`eme~\ref{canon} en
restriction aux feuilles de ${\cal F}_0$.}
\end{rem}

\subsection{D\'emonstration des th\'eor\`emes 3 et 4}    \label{th3}

\hspace{5mm}
Soit $\Lambda_0=({\cal F}_0,\sigma_0)$ une structure de Poisson
minimale sur une vari\'et\'e $M_0$. Pour
d\'emontrer le th\'eor\`eme~\ref{preW},
on suppose tout d'abord que ${\cal P}_0$ est un sous-groupe discret de {\bf R}.
D'apr\`es le th\'eor\`eme~\ref{Poisson}, il existe une pr\'equantification
\( q : (E,\theta) \rightarrow (M,\omega) \)
de l'int\'egration de Poisson $(M,\Lambda)$ qui est triviale en restriction
\`a $M_0$.

D'autre part, le groupo\"{\i}de symplectique
$$ (\Gamma,\eta) = (\nu^\ast{\cal F}, p^\ast \omega - dL)
\stackrel{p}{\longrightarrow} (M,\Lambda) \sgroup{\beta_0}{\alpha_0}
(M_0,\Lambda_0) $$
r\'ealise l'int\'egration symplectique de $(M_0,\Lambda_0)$.
La projection $p$ induit un $S^1$-fibr\'e
principal
\( \widehat{q} : \widehat{E} \rightarrow \Gamma \)
qui est encore trivial en restriction \`a $M_0$ et l'on a le diagramme suivant:

\begin{center}
\begin{picture}(180,45)(0,0)
\put(70,38){\vector(1,0){40}}
\put(70,5){\vector(1,0){40}}
\put(56,29){\vector(0,-1){15}}
\put(123,29){\vector(0,-1){15}}
\put(53,35){$\widehat{E}$}
\put(54,2){$\Gamma$}
\put(90,43){$\widehat{p}$}
\put(90,10){$p$}
\put(119,2){$M$}
\put(120,35){$E$}
\put(131,20){$q$}
\put(44,20){$\widehat{q}$}
\end{picture}
\end{center}

\noindent
La courbure de la connexion relev\'ee $\widehat{p}^\ast \theta$
est \'egale \`a $p^\ast\omega$. Il s'ensuit que la courbure de la connexion
$\widehat{\theta} = \widehat{p}^\ast \theta - \widehat{q}^\ast L$
est \'egale \`a la forme symplectique
$\eta=p^\ast\omega - dL$. Puisque la restriction de $L$ \`a $M_0$ est nulle,
$\widehat{\theta}$ induit la connexion canonique sur le fibr\'e trivial
au-dessus de $M_0$. Bref,
$ \widehat{q} : (\widehat{E},\widehat{\theta}) \rightarrow (\Gamma,\eta)$
est une pr\'equantification triviale en restriction \`a $M_0$.
D'apr\`es le th\'eor\`eme~\ref{canon}, l'espace total
$\widehat{E}$ est muni d'une structure de groupo\"{\i}de de Lie qui en fait
une extension de $\Gamma$ par $S^1$.
\vspace{1ex}

R\'eciproquement, si $(\Gamma,\eta)$ est pr\'equantifiable, le groupe des
p\'eriodes $Per(\eta)$ est un sous-groupe discret de {\bf R}.
Si $\widehat{z}$ est un 2-cycle de $\Gamma$ et $z$ est le 2-cycle projet\'e de
$M$, alors la p\'eriode

\[ \int_{\textstyle \widehat{z}} \eta =
\int_{\textstyle \widehat{z}} p^\ast \omega - dL =
\int_{\textstyle z} \omega \]

\noindent
Puisque
\( \Gamma = \nu^\ast{\cal F} \)
se r\'etracte par d\'eformation sur $M$, les groupes
$Per(\eta)$ et $Per(\omega)$ sont \'egaux. D'autre part,
$\cal P$ est un sous-groupe de $Per(\omega)$, car $\omega$
repr\'esente $\sigma$. Il s'ensuit que $\cal P$ est un sous-groupe discret de
{\bf R}.
Enfin, le groupe ${\cal P}_0$ co\"{\i}ncide avec le groupe $\cal P$
d'apr\`es le lemme~\ref{comp},
ce qui ach\`eve la preuve du th\'eor\`eme~\ref{preW}.
\vspace{2ex}

Pour d\'emontrer le th\'eor\`eme~\ref{tot}, il suffit de remarquer que le
r\'epresentant $\omega$ de $\sigma$ et la forme symplectique
$\eta = p^\ast\omega - dL$ sont exactes pour une vari\'et\'e de Poisson
totalement asph\'erique. Dans ce cas, l'int\'egration symplectique
$(\Gamma,\eta)$ est
pr\'equantifi\'ee par un fibr\'e trivial en cercles.

\subsection{Exemples}

1) Un {\em syst\`eme m\'ecanique} est la donn\'ee d'une vari\'et\'e
symplectique $(S,\omega_S)$ et d'une fonction diff\'erentiable
$H : S \times \mbox{\bf R} \to \mbox{\bf R}$ appel\'ee le
{\em hamiltonien du syst\`eme}. D'apr\`es \cite{A}, on lui associe une {\em
structure cosymplectique}  sur $M = S \times \mbox{\bf R}$ d\'efinie par la
1-forme ferm\'ee $\theta=dt$ et la 2-forme ferm\'ee
$\omega = \omega_S + dH \wedge dt$. C'est l'{\em espace d'\'evolution}
\cite{So} du syst\`eme
m\'ecanique. Pour que celui-ci soit pr\'equantifiable, il faut et il suffit que
le groupe des p\'eriodes de  $\omega$ (qui co\"{\i}ncide avec celui de
$\omega_S$) soit discret.
\vspace{2ex}

\noindent
2) Soit $M$ une vari\'et\'e cosymplectique munie d'une 1-forme
ferm\'ee $\theta$ et d'une 2-forme ferm\'ee $\omega$.
Le feuilletage $\cal F$ d\'efini par l'\'equation $\theta=0$ et
la forme feuillet\'ee $\sigma$ repr\'esent\'ee par $\omega$ d\'efinissent une
structure de Poisson $\Lambda$ sur $M$.
Si le groupe des p\'eriodes $Per(\omega)$
est discret, la vari\'et\'e cosymplectique est pr\'equantifiable.
Pour qu'elle soit pr\'equantifiable au sens de Weinstein, il faut et il suffit
que le groupe $\cal P$ des p\'eriodes sph\'eriques de $\sigma$ soit discret.
Si l'on suppose que le {\em champ de Reeb} $R$
(d\'efini par $i_{\textstyle R} \theta = 1$ et $i_{\textstyle R} \omega = 0$)
est complet, le passage au rev\^etement universel permet de  v\'erifier que
$\cal P$ est \'egal au groupe des p\'eriodes sph\'eriques de $\omega$.

On se propose d'exhiber un exemple de vari\'et\'e cosymplectique
qui n'est pas pr\'equantifiable, mais qui l'est au sens de Weinstein.
\vspace{2ex}

\noindent
3) Soit $\rho$ un nombre irrationel. La forme symplectique
$ dQ_1 \wedge dQ_2 + \rho dq_1 \wedge dq_2 $
sur
$(\mbox{\bf R}^4;Q_1,Q_2,q_1.q_2)$
passe au quotient en une forme symplectique $\omega_S$ sur
$S = T^2 \times T^2$. On consid\`ere la structure cosymplectique sur
$M = S \times \mbox{\bf R}$ d\'etermin\'ee par la 1-forme $\theta=dt$ et la
2-forme
$$ \omega = \omega_S + dH \wedge dt =
dQ_1 \wedge dQ_2 + \rho dq_1 \wedge dq_2 + dH \wedge dt $$
o\`u l'on identifie abusivement la 2-forme $\omega_S$ \`a la 2-forme
relev\'ee sur $\mbox{\bf R}^4$.
Le groupe des p\'eriodes $Per(\omega)$ est partout dense dans {\bf R} et donc
la vari\'et\'e cosymplectique ne poss\`ede pas de pr\'equantification.
\vspace{2ex}

D'autre part, soit $\cal F$ le feuilletage horizontal d\'efini par l'\'equation
$\theta=0$. Le groupo\"{\i}de d'homotopie $\Pi_1({\cal F})$ est
isomorphe au groupo\"{\i}de (voir \cite{We})
$$ T^\ast(T^2) \times T^\ast(T^2) \times \mbox{\bf R}
\sgroup{\beta_0}{\alpha_0} T^2 \times T^2 \times \mbox{\bf R} = M $$
Les projections source et but sont donn\'ees par \vspace{1ex}
$$ \alpha_0(P_1,P_2,Q_1,Q_2,p_1,p_2,q_1,q_2,t) =
(Q_1 + \frac{1}{2}P_2, Q_2 - \frac{1}{2}P_1,
q_1 + \frac{1}{2\rho}p_2, q_2 - \frac{1}{2\rho}p_1,t) \vspace{1ex} $$
$$ \beta_0(P_1,P_2,Q_1,Q_2,p_1,p_2,q_1,q_2,t) =
(Q_1 - \frac{1}{2}P_2, Q_2 + \frac{1}{2}P_1,
q_1 - \frac{1}{2\rho}p_2, q_2 + \frac{1}{2\rho}p_1,t) \vspace{2ex} $$
et le produit de deux \'el\'ements composables
$$
(P_1^\prime,P_2^\prime,Q_1^\prime,Q_2^\prime,
p_1^\prime,p_2^\prime,q_1^\prime,q_2^\prime,t^\prime)
\;\;\;\; \mbox{et} \;\;\;\;
(P_1,P_2,Q_1,Q_2,p_1,p_2,q_1,q_2,t) $$
est \'egal \`a \vspace{1ex}
$$ (P_1^\prime + P_1,P_2^\prime + P_2, Q_1 - \frac{1}{2}P_2^\prime,
Q_2 + \frac{1}{2}P_1^\prime,
p_1^\prime + p_1, p_2^\prime + p_2, q_1 - \frac{1}{2\rho}p_2^\prime,
q_2 + \frac{1}{2\rho}p_1^\prime,t^\prime + t) \vspace{2ex} $$
La 2-forme ferm\'ee
$ \widetilde{\omega} = \alpha_0^\ast \omega - \beta_0^\ast \omega $
s'\'ecrit \vspace{1ex}
$$\begin{array}{lcl}
\widetilde{\omega} & = &
dP_1 \wedge dQ_1 + dP_2 \wedge dQ_2
+ dp_1 \wedge dq_1 + dp_2 \wedge dq_2
+ \alpha_0^\ast(dH \wedge dt) - \beta_0^\ast(dH \wedge dt) \\
 & &  \\
 & = &
d(P_1 dQ_1 + P_2 dQ_2 + p_1 dq_2 + p_2 dq_2 + \alpha_0^\ast(H dt) -
\beta_0^\ast(H dt)) \\
 & &  \\
 & = &
d(L + l + \alpha_0^\ast(H dt) - \beta_0^\ast(H dt))
\end{array} \vspace{2ex} $$
La 1-forme ferm\'ee $\widetilde{\theta}=dt$ et la 2-forme ferm\'ee
$\widetilde{\omega}$
d\'efinissent un structure de groupo\"{\i}de cosymplectique sur
$ T^\ast(T^2) \times T^\ast(T^2) \times \mbox{\bf R} $.
Le groupo\"{\i}de
$$ \Gamma = T^\ast(T^2) \times T^\ast(T^2) \times T^\ast(\mbox{\bf R})
\stackrel{p}{\longrightarrow}
T^\ast(T^2) \times T^\ast(T^2) \times \mbox{\bf R}
\sgroup{\beta_0}{\alpha_0} T^2 \times T^2 \times \mbox{\bf R} = M $$
muni de la forme symplectique
$$ \eta = p^\ast \widetilde{\omega} - dr \wedge dt =
d(L + l + \alpha^\ast(H dt) - \beta^\ast(H dt)) - rdt) = d \lambda $$
r\'ealise l'int\'egration symplectique de la vari\'et\'e cosymplectique
$M$, o\`u l'on note $\alpha$ et $\beta$ les projections source et but.
La vari\'et\'e cosymplectique $M$ est pr\'equantifi\'ee au
sens de Weinstein par le fibr\'e trivial
$ \widehat{q} : \Gamma \times S^1 \rightarrow \Gamma$
muni de la connexion
$dz + \widehat{q}^\ast \lambda$, o\`u $dz$ est la connexion canonique.
\vspace{3ex}

\noindent
{\bf Remerciements} C'est gr\^ace \`a de nombreuses discussions avec
Gilbert Hector que ce travail a pu \^etre r\'ealis\'e. Je lui suis
tr\`es reconnaissant.
Je veux aussi remercier Gijs Tuynman pour ses remarques.

\end{document}